\documentclass[11pt,notitlepage]{article}
\usepackage[margin=1in]{geometry}
\usepackage[parfill]{parskip}   
\usepackage[nodisplayskipstretch]{setspace}
\usepackage{enumerate}
\usepackage{bm}
\usepackage{amssymb}
\usepackage{amsmath}

\usepackage{amsthm}
\usepackage{graphicx}
\usepackage{subfigure}
\usepackage{makeidx}
\usepackage{multicol}
\usepackage{comment}
\usepackage{color}
\usepackage{chngcntr}
\definecolor{Gray}{gray}{0.9}
\newtheorem{theorem}{Theorem}
\newtheorem*{mproof}{Proof}
\newtheorem{lemma}{Lemma}
\newtheorem{proposition}{Proposition}
\newtheorem{corollary}{Corollary}
\newtheorem{remark}{Remark}



\date{}

\begin{document}

\vspace{-10pt}
\title{\large\bf  On The Number Of Unlabeled Bipartite Graphs}

\vspace{-30pt}
\author{\ Abdullah Atmaca\footnote{Department of Computer Engineering, Bilkent University, Ankara, Turkey.}   \,
  and A. Yavuz Oru\c{c}\footnote{Department of Electrical and Computer Engineering, University of Maryland, College Park, MD 20742.}
}

\maketitle
\nopagebreak[4]
\setcounter{page}{1}

\vspace{-25pt}
\begin{center}
{\bf Abstract}
\end{center}

\vspace{-10pt}\noindent
This paper solves a problem that was stated  by M. A. Harrison in 1973~\cite{harrison1973number}. This problem, that has remained open since then is concerned with counting  equivalence classes of $n\times r$ binary matrices under row and column permutations.  Let  $I$ and $O$ denote two sets of vertices, where $I\cap O =\Phi$,  $|I| = n$, $|O| = r$,  and  $B_u(n,r)$ denote the set of unlabeled graphs whose edges connect vertices in $I$ and $O$.  Harrison established that the number of equivalence classes of $n\times r$ binary matrices is equal to the number of unlabeled graphs in $B_u(n,r).$ He  also computed the number of such matrices (hence such graphs) for small values of $n$ and $r$ without providing an asymptotic formula $|B_u(n,r)|.$ 
Here, such an asymptotic formula is provided  by proving the following two-sided equality using Polya's Counting Theorem.

\vspace{-18pt}
\begin{equation}
\label{mainResult}
\displaystyle  \frac{\binom{r+2^{n}-1}{r}}{n!} \le |B_u(n,r)| \le  2\frac{\binom{r+2^{n}-1}{r}}{n!} 
\end{equation}

\vspace{-20pt}\noindent
where $n < r$. 

\vspace{-20pt}

\section{Introduction}
\label{problemStatement}

\vspace{-17pt}\noindent
The counting problem that is considered in this paper
has been investigated in connection with the enumeration of  unlabeled bipartite graphs and binary matrices\cite{harrison1973number}. Let $(I,O,E)$ denote a  graph with two disjoint  sets of vertices, $I$, called {\em left vertices} and a set of vertices, $O$, called {\em right vertices}, where  each edge in $E$ connects a left vertex with a right vertex. We let $n = |I|$, $r = |O|$, and refer to such a graph as an $(n,r)$-bipartite graph. Let  $G_1 = (I,O,E_1)$ and $G_2 = (I,O,E_2)$ be two $(n,r)$-bipartite graphs, and $\alpha: I\rightarrow I$  and $\beta: O\rightarrow O$ be both bijections. The pair $(\alpha,\beta)$ is  an isomorphism between $G_1$ and $G_2$ provided that $((\alpha(v_1), \beta(v_2))\in E_2$ if and only if $(v_1,v_2)\in E_1$,  $\forall v_1\in I, \forall v_2\in O$.  It is easy to establish that this mapping induces an equivalence relation, and partitions the set of $2^{nr}$  $(n,r)$-bipartite graphs into equivalence classes. This equivalence relation captures the fact that  the vertices in $I$ and $O$ are unlabeled, and  so each class of $(n,r)$-bipartite graphs can be represented by any one of the graphs in that class without identifying the vertices in $I$ and $O$. Let $B_u(n,r)$ denote any set of $(n,r)$-bipartite graphs that  contains  exactly one such graph from each of the equivalence classes of $(n,r)$-bipartite graphs induced by the isomorphism we defined. It is easy to see that determining $|B_u(n,r)|$ amounts to an enumeration of non-isomorphic  $(n,r)$-bipartite graphs that will henceforth be referred to as unlabeled $(n,r)$-bipartite graphs. 

\vspace{-10pt}
In \cite{harrison1973number}, Harrison used P\'{o}lya's counting theorem  to obtain an expression to compute the number of non-equivalent $n\times r$ binary matrices. This expression contains a nested sum, in which  one sum is carried over all partitions of $n$ while the other is carried over all partitions of $r$, where the argument of the nested sum involves  factorial, exponentiation and greatest common divisor (gcd) computations. He further established that this formula also enumerates the number of  unlabeled $(n,r)$-bipartite graphs, i.e., $|B_u(n,r)|$. 
A number of results indirectly related to Harrison's work and our result appeared in the literature~\cite{harary1958number,harary1963enumeration,hanlon1979enumeration}. In particular, the set $B_u(n,r)$ in our work coincides with the set of bicolored graphs described in Section 2 in~\cite{harary1958number}. Whereas \cite{harary1958number} provides a counting polynomial for the number of  bicolored graphs, we focus on the asymptotic behavior of $|B_u(n,r)|$ in this paper.  Counting polynomials for other families of  bipartite graphs were also reported in  \cite{harary1963enumeration}. Likewise, \cite{hanlon1979enumeration} provides generating functions for related bipartite graph counting problems without an asymptotic analysis as provided in this paper. It should also be mentioned that some results on asymptotic enumeration of certain families of bipartite graphs (binary matrices) have been reported (see for example, \cite{canfield2005asymptotic,greenhill2006asymptotic,canfield2008asymptotic,barvinok2010number}). To the best of our knowledge, our work provides the first asymptotic enumeration of unlabelled bipartite graphs.

\vspace{-8pt}
Let $S_n$ denote the symmetric group of permutations of degree $n$ acting on set $N =\{1,2,\cdots, n\}$. Suppose that the $n!$ permutations in $S_n$ are indexed by $1,2,\cdots,n!$ in some arbitrary, but fixed manner.  The cycle index polynomial of $S_n$ is defined as follows(\cite{hararyEnumeration},see p.35, Eqn. 2.2.1):
\vspace{-5pt}
\begin{equation}
Z_{S_n} (x_1,x_2,\cdots, x_n) = \frac{1}{n!}\sum^{n!}_{m=1} \prod_{k=1}^n x_k^{p_{m,k}}
\end{equation}

\vspace{-15pt}\noindent
where $p_{m,k}$ denotes the number of cycles of length $k$ in the disjoint cycle representation of the  $m^{\rm th}$ permutation in $S_n$, and $\sum_{k=1}^n kp_{m,k} = n, \forall m = 1,2,\cdots, n!.$

\vspace{1pt}
Let $S_{n} \times S_{r}$ denote the direct product of symmetric groups $S_{n}$ and $S_{r}$ acting on $N =\{1,2,\cdots, n\}$ and $R =\{1,2,\cdots, r\},$ respectively, where $n$ and $r$ are positive integers such that $n < r.$ It can be inferred from Harrison (\cite{HarrisonChapter4},Lemma 4.1 and Theorem 4.2) that the cycle index polynomial of $S_{n} \times S_{r}$ is given by \cite{HarrisonChapter4}

\vspace{-5pt}
\begin{equation}
\label{HarrisonEquation1}
Z_{S_{n} \times S_{r}}(x_1,x_2,\cdots, x_{nr}) = Z_{S_{n}}(x_1,x_2,\cdots, x_{n}) \boxtimes Z_{S_{r}}(x_1,x_2,\cdots, x_{r}),
\end{equation}

\vspace{-15pt}\noindent 
where $\boxtimes$ is a particular polynomial multiplication that distributes over ordinary addition, and in which the multiplication  $X_m \bigodot X_t$  of two product terms\footnote{Note that we will not display the zero powers of $x_1, x_2,\cdots$   
 in a cycle index polynomial. We will use the same convention for all other cycle index polynomials throughout the paper.}, $X_{m}=x_{1}^{p_{m,1}}x_{2}^{p_{m,2}} \cdots x_{n}^{p_{m,n}}$ and $X_{t}=x_{1}^{q_{t,1}}x_{2}^{q_{t,2}} \cdots x_{r}^{q_{t,r}}$  in $Z_{S_{n}}$ and $Z_{S_{r}}$, respectively,  is defined as\footnote{The lcm($a$,$b$) and gcd($a$,$b$) denote least common multiple and greatest common divisor of $a$ and $b$.}

\vspace{-25pt}
\begin{equation}
X_{m}\bigodot X_{t} 
= \prod_{k=1}^{n} \prod_{j=1}^{r} x_{{\rm lcm}(k,j) }^{p_{m,k}q_{t,j} {\rm gcd}(k,j)}.
\end{equation}

\vspace{-10pt}
Harrison further proved that \cite{harrison1973number}:

\vspace{-20pt}
\begin{equation}
\label{HarrisonGraphCount}
\vert B_{u}(n,r)\vert = Z_{S_{n} \times S_{r}}(\underbrace{2,2,.. , 2}_{nr})
\end{equation}

\vspace{-25pt}\noindent
when\footnote{As noted in \cite{harrison1973number}, $n = r$ case involves a different cycle index polynomial.  Bounding $|B_u(n,n)|$ will be considered separately at the end of the paper.} $n \neq r$.

\vspace{-10pt}
We need one more fact that can be found in Harary (\cite{hararyEnumeration}, p.36) in order to compute the stated lower and upper bound in (\ref{mainResult}):

\vspace{-25pt}
\begin{equation}
\label{cycleIndexRecurrence}
Z_{S_{r}}(x_{1},x_{2},\ldots\ldots,x_{r}) =
	\frac{1}{r} \sum_{i=1}^{r}x_{i} Z_{S_{r-i}}(x_{1},x_{2},\ldots\ldots,x_{r-i})
\end{equation}

\vspace{-10pt}\noindent
where $Z_{S_{0}}() = 1$.

\section{The Lower Bound for $\vert B_{u}(n,r)\vert$}

\vspace{-15pt}

From (\ref{HarrisonEquation1}) and (\ref{HarrisonGraphCount}) we know that

\vspace{-20pt}\noindent
\begin{align}
& \vert B_{u}(n,r)\vert & =& Z_{S_{n} \times S_{r}}(2,2,\ldots,2), \\
&& =& [Z_{S_{n}}(x_1,x_2,\cdots, x_{n}) \boxtimes Z_{S_{r}}(x_1,x_2,\cdots, x_{r}) ](2,2,\ldots,2).
\end{align}

\vspace{-20pt}\noindent
One of the terms in $Z_{S_{n}}(x_1,x_2,\cdots, x_{n})$ is $\frac{1}{n!}(x_{1}^{n})$ and it is associated with the identity permutation in $S_{n}$. Using this fact, we find

\vspace{-20pt}\noindent
\begin{align}
& \vert B_{u}(n,r)\vert & =& Z_{S_{n} \times S_{r}}(2,2,\ldots,2), \\
&& =& [Z_{S_{n}}(x_1,x_2,\cdots, x_{n}) \boxtimes Z_{S_{r}}(x_1,x_2,\cdots, x_{r}) ](2,2,\ldots,2), \\
&& =& \left [ \left (\frac{1}{n!}(x_{1}^{n}+\ldots) \right ) \boxtimes Z_{S_{r}}(x_1,x_2,\cdots, x_{r})\right ] (2,2,\ldots,2), \\
&& =& \left[ \left (\frac{1}{n!}x_{1}^{n} \right ) \boxtimes Z_{S_{r}}(x_{1},x_{2},\ldots,x_{r}) \right](2,2,\ldots,2) + \ldots,  \\
&& =& \frac{1}{n!} \Bigg\lbrace \Big [x_{1}^{n}\boxtimes \frac{1}{r!} \sum_{t=1}^{r!} \prod_{j=1}^{r}
x_{j}^{q_{t,j}}\Big ] (2,2,...,2) \Bigg\rbrace + \ldots,\\ 
&&=& \frac{1}{n!} \Bigg\lbrace \Big [\frac{1}{r!} \sum_{t=1}^{r!} x_{1}^{n} \bigodot \prod_{j=1}^{r} x_{j}^{q_{t,j}}\Big ] (2,2,...,2) \Bigg\rbrace + \ldots, \\
&&=& \frac{1}{n!}  \Bigg\lbrace \Big [\frac{1}{r!} \sum_{t=1}^{r!} \prod_{j=1}^{r} x_{{\rm lcm}(1,j) }^{n q_{t,j} {\rm gcd}(1,j)  }\Big ] (2,2,...,2) \Bigg\rbrace + \ldots, \\
&& =& \frac{1}{n!} \Bigg\lbrace \Big [\frac{1}{r!} \sum_{t=1}^{r!} \prod_{j=1}^{r}  x_{j}^{n q_{t,j}}
\Big ] (2,2,...,2) \Bigg\rbrace + \ldots,\\
&& =& \frac{1}{n!} \Bigg\lbrace \frac{1}{r!} \sum_{t=1}^{r!} \prod_{j=1}^{r}  2^{n q_{t,j}} \Bigg\rbrace + \ldots, \\
&& =& \frac{1}{n!} \Bigg\lbrace \frac{1}{r!} \sum_{t=1}^{r!} \prod_{j=1}^{r}  (2^{n})^{q_{t,j}} \Bigg\rbrace + \ldots,\\
&&  = & \frac{1}{n!} \Bigg\lbrace Z_{S_{r}}(2^{n},2^{n},\ldots,2^{n}) \Bigg\rbrace + \ldots.
\end{align}

\vspace{-20pt}\noindent
This proves

\vspace{-30pt}\noindent
\begin{align}
\label{ZSrInequality}
& \vert B_{u}(n,r)\vert \geq \frac{1}{n!} Z_{S_{r}}(2^{n},2^{n},\ldots,2^{n}). \qed
\end{align}

\vspace{-20pt}\noindent
\begin{proposition}
\label{lowerBoundProp1}
$$Z_{S_{r}}(2^{n},2^{n},\ldots,2^{n}) = \binom{r+2^{n}-1}{r}$$
\end{proposition}

\vspace{-20pt}\noindent
\begin{mproof}
{\rm
Using (\ref{cycleIndexRecurrence}), we have

\vspace{-28pt}\noindent
\begin{align}
\label{lowerBoundProp1Eq1}
&&  r Z_{S_{r}}(2^{n},2^{n},\ldots,2^{n}) = \sum_{i=1}^{r}2^{n} Z_{S_{r-i}}(2^{n},2^{n},\ldots,2^{n}),
\end{align}

\vspace{-20pt}\noindent
and
\vspace{-20pt}\noindent
\begin{align}
\label{lowerBoundProp1Eq2}
&&  (r-1) Z_{S_{r-1}}(2^{n},2^{n},\ldots,2^{n}) = \sum_{i=1}^{r-1}2^{n} Z_{S_{r-1-i}}(2^{n},2^{n},\ldots,2^{n}).
\end{align}

\vspace{-20pt}\noindent
Subtracting the second equation from the first one 
gives 

\vspace{-20pt}\noindent
\begin{align}
\label{lowerBoundProp1Rec1}
r Z_{S_{r}}(2^{n},2^{n},\ldots,2^{n}) - (r-1) Z_{S_{r-1}}(2^{n},2^{n},\ldots,2^{n}) &= 2^{n} Z_{S_{r-1}}(2^{n},2^{n},\ldots,2^{n}),  \\
r Z_{S_{r}}(2^{n},2^{n},\ldots,2^{n}) &= (r+2^{n}-1) Z_{S_{r-1}}(2^{n},2^{n},\ldots,2^{n}), \\
Z_{S_{r}}(2^{n},2^{n},\ldots,2^{n}) &= (\frac{r+2^{n}-1}{r}) Z_{S_{r-1}}(2^{n},2^{n},\ldots,2^{n}).
\end{align}

\vspace{-23pt}\noindent
Expanding the last equation inductively, we obtain
%Now we will rewrite $Z_{S_{r-1}}(2^{n},2^{n},\ldots,2^{n})$ in Eq.\ref{lowerBoundProp1Rec1} as following:

\vspace{-18pt}\noindent
\begin{align}
Z_{S_{r}}(2^{n},2^{n},\ldots,2^{n}) &= (\frac{r+2^{n}-1}{r})(\frac{r+2^{n}-2}{r-1}) Z_{S_{r-2}}(2^{n},2^{n},\ldots,2^{n}),\\
Z_{S_{r}}(2^{n},2^{n},\ldots,2^{n}) &= (\frac{r+2^{n}-1}{r})(\frac{r+2^{n}-2}{r-1})(\frac{r+2^{n}-3}{r-2}) Z_{S_{r-3}}(2^{n},2^{n},\ldots,2^{n}),\\
Z_{S_{r}}(2^{n},2^{n},\ldots,2^{n}) &= (\frac{r+2^{n}-1}{r})(\frac{r+2^{n}-2}{r-1})(\frac{r+2^{n}-3}{r-2}) \ldots (\frac{2^{n}}{1}) Z_{S_{0}}().
\end{align}

\vspace{-30pt}\noindent

\noindent
Noting that $Z_{S_{0}}() = 1,$ and combining the product terms together, we obtain

\vspace{-20pt}\noindent
\begin{eqnarray}
Z_{S_{r}}(2^{n},2^{n},\ldots,2^{n}) = \binom{r+2^{n}-1}{r}.\qed
\end{eqnarray}
}
\end{mproof}

\vspace{-20pt}\noindent
Combining Proposition \ref{lowerBoundProp1} with (\ref{ZSrInequality}) proves the lower bound:

\vspace{5pt}
\begin{theorem}
$\,$\\
\vspace{-25pt}
\label{theorem1-1}
\begin{align}
\vert B_{u}(n,r)\vert & \geq \frac{1}{n!} Z_{S_{r}}(2^{n},2^{n},\ldots,2^{n})  \geq \frac{ \binom{r+2^{n}-1}{r}}{n!}.\qed
\end{align}
\end{theorem}

\vspace{-25pt}
\section{An Upper Bound for $\vert B_{u}(n,r)\vert$}

\vspace{-15pt}
We first note that $|B_u(1,r)| = r +1= {r+2^1 -1\choose r}/1! \le 2 {r+2^1 -1\choose r}/1!$. Hence the upper bound that is claimed in the abstract holds for $n=1.$ Proving that it also holds for $n\ge 2$ requires a more careful analysis of the terms in 

\vspace{-20pt}
\begin{equation}
\label{HarrisonEquation2}
Z_{S_{n}}(x_1,x_2,\cdots, x_{n}) \boxtimes Z_{S_{r}}(x_1,x_2,\cdots, x_{r}).
\end{equation}
We first express $Z_{S_{n}}(x_1,x_2,\cdots, x_{n}) $  as 

\vspace{-30pt}
\begin{align}
& Z_{S_{n}}(x_{1},x_{2},\ldots,x_{n}) = Z_{S_{n}}[1] + Z_{S_{n}}[2] + \ldots + Z_{S_{n}}[n!],
\end{align}

\vspace{-25pt}\noindent
where

\vspace{-35pt}
\begin{align}
& Z_{S_{n}}[1] =  \frac{1}{n!} x_{1}^{n} \qquad\quad 
\\
& Z_{S_{n}}[2] =  \frac{1}{n!} x_{1}^{n-2}x_{2} \quad 
\end{align}

\vspace{-20pt}\noindent
The first term is associated with  the identity permutation and the second term is associated with any one of the permutations in which all but two of the elements in $N={1,2,\cdots, n}$ are fixed to themselves.
The remaining $Z_{S_{n}}[i] =\frac{1}{n!} \prod_{k=1}^{n} x_{k}^{p_{i,k}}, 3\le i\le n!$ terms represent all the other product terms in the cycle index polynomial of $S_{n}$ with no particular association with the permutations in $S_{n}.$ Similarly, we set $Z_{S_{r}}(x_{1},x_2,\ldots,x_{r}) = \frac{1}{r!} \sum_{t=1}^{r!}  \prod_{j=1}^{r} x_{j}^{q_{t,j}}$ without identifying the actual product terms with any particular permutation in $S_{r}$.

\vspace{-10pt}\noindent\noindent The following equations obviously hold as the sum of the lengths of all the cycles in any cycle disjoint representation of a permutation in $S_{n}$ and $S_{r}$ must be $n$ and $r,$ respectively.

\vspace{-20pt}\noindent
\begin{align}
\label{GroupTheoryEquations1}
& \sum_{k=1}^{n} k p_{i,k} = n, 1 \leq i \leq n!, \\
& \sum_{j=1}^{r} j q_{t,j} = r, 1 \leq t \leq r!
\end{align}

\vspace{-20pt}\noindent Now we can proceed with the computation of the upper bound for $\vert B_{u}(n,r)\vert$. First, we note that

\vspace{-20pt}\noindent 
\begin{align}
\vert B_{u}(n,r)\vert =& Z_{S_{n} \times S_{r}}(2,2,2,\ldots,2), \\
=& \left [ Z_{S_{n}}(x_{1},x_{2},\ldots,x_{n}) \boxtimes Z_{S_{r}}(x_1,x_2,\cdots, x_{r})\right ] (2,2,\ldots,2),\\
=& \left [ \left(Z_{S_{n}}[1] + Z_{S_{n}}[2] + \ldots + Z_{S_{n}}[n!]\right) \boxtimes Z_{S_{r}}(x_{1},x_{2},\ldots,x_{r}) \right ] (2,2,\ldots,2),\\
=& \left [Z_{S_{n}}[1] \boxtimes Z_{S_{r}}(x_{1},x_2,\ldots,x_{r})\right ] (2,2,\ldots,2) 
 + \left [Z_{S_{n}}[2] \boxtimes
Z_{S_{r}}(x_{1},x_2,\ldots,x_{r})\right ] (2,2,\ldots,2)\nonumber\\
& + \ldots + \left [Z_{S_{n}}[n!] \boxtimes Z_{S_{r}}(x_{1},x_2,\ldots,x_{r})\right ] (2,2,\ldots,2).\label{firstTerm}
\end{align}

\vspace{-15pt}\noindent 
The first term in (\ref{firstTerm}) is directly computed from Proposition~\ref{lowerBoundProp1}. Thus, it suffices to upper bound each of the remaining terms in (\ref{firstTerm}) to upper bound $|B_{u}(n,r)|.$ This will be established by proving $ \left [Z_{S_{n}}[2] \boxtimes
Z_{S_{r}}(x_{1},x_2,\ldots,x_{r})\right ] (2,2,\ldots,2)\ge
 \left [Z_{S_{n}}[i] \boxtimes Z_{S_{r}}(x_{1},x_2,\ldots,x_{r})\right ] (2,2,\ldots,2), \forall i, 3\le i\le n!.$ We first need some preliminary facts.

\vspace{5pt}
\begin{lemma}
\label{lemmaProp5}
{\rm
 For all $i, 1\le i\le n!$, }
 
 \vspace{-20pt}\noindent
\begin{eqnarray}
& [Z_{S_{n}}[i] \boxtimes Z_{S_{r}}(x_{1},x_2,\ldots,x_{r})](2,\ldots,2) = \frac{1}{n!}Z_{S_{r}}(2^{\sum_{k=1}^{n}p_{i,k} {\rm gcd}(k,1)},\ldots,2^{\sum_{k=1}^{n}p_{i,k} {\rm gcd}(k,r)}).
\label{lemmaProp5Eq1}
\end{eqnarray}
\end{lemma}

\vspace{-10pt}
\begin{mproof}

$\,$\\
\vspace{-20pt}
\begin{align}
[Z_{S_{n}}[i] \boxtimes Z_{S_{r}}(x_{1},x_2,\ldots,x_{r})](2,2,\ldots,2) & =\left[\frac{1}{n!} \prod_{k=1}^{n} x_{k}^{p_{i,k}} \boxtimes \left (\frac{1}{r!} \sum_{t=1}^{r!} \prod_{j=1}^{r} x_{j}^{q_{t,j}} \right ) \right](2,2,\ldots,2), \\
& = \left[\frac{1}{n! r!} \sum_{t=1}^{r!} \prod_{k=1}^{n}  x_{k}^{p_{i,k}}  \bigodot   \prod_{j=1}^{r} x_{j}^{q_{t,j}}  \right](2,2,\ldots,2),\\
& = \left[\frac{1}{n! r!} \sum_{t=1}^{r!} \prod_{j=1}^{r} \prod_{k=1}^{n} x_{{\rm lcm}\left ( k,j \right )}^{p_{i,k} q_{t,j} {\rm gcd}(k,j) } \right](2,2,\ldots,2), \\
& = \frac{1}{n! r!} \sum_{t=1}^{r!} \prod_{j=1}^{r} \prod_{k=1}^{n} 2^{p_{i,k} q_{t,j} {\rm gcd}(k,j) },\\
& = \frac{1}{n!} \left[\frac{1}{r!} \sum_{t=1}^{r!} \prod_{j=1}^{r}  (2^{\sum_{k=1}^{n} p_{i,k} {\rm gcd}(k,j) })^{q_{t,j}} \right],  \\
& = \frac{1}{n!}Z_{S_{r}}(2^{\sum_{k=1}^{n}p_{i,k} {\rm gcd}(k,1)},\ldots,2^{\sum_{k=1}^{n}p_{i,k} {\rm gcd}(k,r)}).\qed 
\end{align}
\end{mproof}

\vspace{-20pt}
\begin{corollary}
\label{corollaryProp5}
$\,$\\

\vspace{-45pt}
\begin{eqnarray}
[Z_{S_{n}}[2] \boxtimes Z_{S_{r}}(x_{1},x_2,\ldots,x_{r})](2,\ldots,2)
&= \frac{1}{n!}Z_{S_{r}}(2^{n-1},2^{n},2^{n-1},2^{n},\ldots).
\label{corollaryProp5Eq1}
\end{eqnarray}
\end{corollary}
\vspace{-10pt}
\begin{mproof}
{\rm By definition, $p_{2,1} = n-2, p_{2,2} = 1, p_{2,k} = 0, 3\le k\le n.$ Substituting these into the last equation in Lemma~\ref{lemmaProp5} proves the statement. \qed}
\end{mproof}

\vspace{5pt}
\begin{lemma}
\label{prop5Lemma2}
{\rm

} 
$\,$
\vspace{-10pt}\noindent
\begin{eqnarray}
& \sum_{k=1}^{n} p_{i,k} \leq n-1, \forall i, 2\le i\le n!.
\label{prop5Lemma2Equation}
\end{eqnarray}
\end{lemma}

\vspace{-5pt}
\begin{mproof}

{\rm
Recall from (\ref{GroupTheoryEquations1}) that $\sum_{k=1}^{n}k p_{i,k} = n,$  $\forall i, 1\le i\le n!.$  Hence $\sum_{k=1}^{n}p_{i,k} = n - \sum_{k=1}^{n}(k-1)  p_{i,k},$ and so the maximum value of $\sum_{k=1}^{n}p_{i,k}$ occurs when $\sum_{k=1}^{n}(k-1)  p_{i,k}$ is minimized. Furthermore, at least one of  $ p_{i,k},\forall i, 2\le i\le n!$ must be $\ge 1$ for some $k\ge 2$ since none of the permutations we consider is the identity. Thus, $\sum_{k=1}^{n}(k-1)  p_{i,k}\ge 1$ and the statement follows. \qed
}
\end{mproof}

\vspace{5pt}
\begin{lemma}
\label{prop5Lemma3}
{\rm

If $\sum_{k=1}^{n}p_{i,k} {\rm gcd}(k,\alpha+1) = n$, then $\sum_{k=1}^{n}p_{i,k} {\rm gcd}(k,\alpha) \leq n-1$, $\forall i, 2\le i\le n!$ and for any integer $\alpha \geq 2$.

}
\end{lemma}
\vspace{8pt}
\begin{mproof}
{\rm
\noindent 
If $\sum_{k=1}^{n}p_{i,k} {\rm gcd}(k,\alpha+1) = n$ as stated in the lemma, then we must have $ {\rm gcd}(k,\alpha+1) = k $   where $p_{i,k}\ge 1$, $\forall i, 2\le i\le n!$. Therefore 
$ k \leq \alpha+1 $. Now if $k=\alpha+1$, then trivially ${\rm gcd}(k,\alpha) < k $. On the other hand if
 $k < \alpha+1$, then  $\alpha+1$ must  be a multiple of $k$. Therefore, $\alpha$ can not be a multiple of $k$ for any $k \ge 2$. At this point we find that ${\rm gcd}(k,\alpha) < k,$ $\forall k, 2\le k\le n$. Since as in the previous lemma, none of the permutations we consider is the identity, at least one of  $ p_{i,k},\forall i, 2\le i\le n!$ must be $\ge 1$ for some $k\ge 2$ and so we conclude that $\sum_{k=1}^{n}p_{i,k} {\rm gcd}(k,\alpha) \leq n-1$.
\qed
}
\end{mproof}
\vspace{5pt}
\begin{lemma}
\label{prop5Lemma4}
{\rm
\noindent 
$Z_{S_{r}}(2^{n-1},2^{n},\ldots) \geq Z_{S_{r-1}}(2^{n-1},2^{n},\ldots)$, for $2 \leq n$.
}
\end{lemma}
\vspace{5pt}
\begin{mproof}
{\rm
\noindent Using (\ref{cycleIndexRecurrence}), we get

\vspace{-35pt}
\begin{align}
\label{Lemma4Eqq1}
& r Z_{S_{r}}(2^{n-1},2^{n},\ldots) =  \sum_{\textbf{odd } i}^{r-\beta_{1}}2^{n-1} Z_{S_{r-i}}(2^{n-1},2^{n},\ldots) + \sum_{\textbf{even } i}^{r-\beta_{2}}2^{n} Z_{S_{r-i}}(2^{n-1},2^{n},\ldots),
\end{align}

\vspace{-15pt}\noindent
where $\beta_{1}=1, \beta_{2}=0$ if $r$ is even and $\beta_{1}=0, \beta_{2}=1$ if $r$ is odd. Similarly, for $r-1$,

\vspace{-35pt}\noindent

\begin{align}
\label{Lemma4Eqq2}
& (r-1) Z_{S_{r-1}}(2^{n-1},2^{n},\ldots) =  \sum_{\textbf{odd } i}^{r-1-\beta_{2}}\!\!\!2^{n-1} Z_{S_{r-1-i}}(2^{n-1},2^{n},\ldots)\!+\!\!\! \sum_{\textbf{even } i}^{r-1-\beta_{1}}\!\!\!2^{n} Z_{S_{r-1-i}}(2^{n-1},2^{n},\ldots).
\end{align}

\vspace{-20pt}\noindent
Subtracting \ref{Lemma4Eqq2} from \ref{Lemma4Eqq1} gives
\vspace{-5pt}\noindent
\begin{align}
& r Z_{S_{r}}(2^{n-1},2^{n},\ldots) - (r-1) Z_{S_{r-1}}(2^{n-1},2^{n},\ldots)\nonumber\\
& \quad =  \sum_{\textbf{even } i}^{r-\beta_{2}}2^{n} Z_{S_{r-i}}(2^{n-1},2^{n},\ldots) - \sum_{\textbf{odd } i}^{r-1-\beta_{2}}2^{n-1} Z_{S_{r-1-i}}(2^{n-1},2^{n}\ldots)\nonumber \\
& \quad \quad + \sum_{\textbf{odd } i}^{r-\beta_{1}}2^{n-1} Z_{S_{r-i}}(2^{n-1},2^{n},\ldots) - \sum_{\textbf{even } i}^{r-1-\beta_{1}}2^{n} Z_{S_{r-1-i}}(2^{n-1},2^{n},\ldots), \\
& r Z_{S_{r}}(2^{n-1},2^{n},\ldots) - (r-1) Z_{S_{r-1}}(2^{n-1},2^{n},\ldots)\nonumber\\
& \quad = \!\!\!\! \sum_{\textbf{even } i}^{r-\beta_{2}}\!\!\!2^{n-1} Z_{S_{r-i}}(2^{n-1},2^{n},\ldots)  + 2^{n-1} Z_{S_{r-1}}(2^{n-1},2^{n},\ldots) - \!\!\!\!\sum_{\textbf{even } i}^{r-1-\beta_{1}}\!\!\!\!2^{n-1} Z_{S_{r-1-i}}(2^{n-1},2^{n},\ldots),
\end{align}
\begin{align}
& r Z_{S_{r}}(2^{n-1},2^{n},\ldots) = (r-1 + 2^{n-1}) Z_{S_{r-1}}(2^{n-1},2^{n}\ldots)\nonumber \\
& \quad  + 2^{n-1} \left( \sum_{\textbf{even } i}^{r-\beta_{2}} Z_{S_{r-i}}(2^{n-1},2^{n},\ldots) -\! \sum_{\textbf{even } i}^{r-1-\beta_{1}} Z_{S_{r-1-i}}(2^{n-1},2^{n},\ldots) \right). \label{inductionEquationLemma4}
\end{align}

\vspace{-15pt}\noindent
We now prove the lemma by induction on $r$.

\vspace{-10pt}\noindent
\textbf{Basis $r=1.$} By (\ref{cycleIndexRecurrence}), $Z_{S_{1}}(2^{n-1}) = 2^{n-1} Z_{S_{0}}() = 2^{n-1}.$ So we have $Z_{S_{1}}(2^{n-1}) = 2^{n-1} \geq Z_{S_{0}}() = 1$ for $2 \leq n$.

\vspace{-10pt}\noindent
\textbf{Induction Step.} Suppose that the lemma holds from $1$ to $r-1$. That is, $Z_{S_{r-i}}  -Z_{S_{r-i-1}} \ge 0,  1\le i\le r-1.$ Now if $r$ is even then the difference of the two sums in~(\ref{inductionEquationLemma4}) becomes $(Z_{S_{r-2}} -Z_{S_{r-3}}) + (Z_{S_{r-4}} -Z_{S_{r-5}})\ldots + (Z_{S_2} -Z_{S_1}) + Z_{S_0}$, which is clearly $\ge 0$ by the induction hypothesis. Therefore, 
\vspace{-10pt}\noindent
\begin{align}
& r Z_{S_{r}}(2^{n-1},2^{n},\ldots) \geq (r-1 + 2^{n-1}) Z_{S_{r-1}}(2^{n-1},2^{n},\ldots),\\
& Z_{S_{r}}(2^{n-1},2^{n},\ldots) \geq Z_{S_{r-1}}(2^{n-1},2^{n},\ldots), n\ge 2.
\end{align}
On the other hand, if $r$ is odd then the difference of the two sums in the same equation becomes $(Z_{S_{r-2}} -Z_{S_{r-3}}) + (Z_{S_{r-4}} -Z_{S_{r-5}})\ldots + (Z_{S_2} -Z_{S_1}) + (Z_{S_1} - Z_{S_0})$, which is again $\ge 0,$ and the statement follows in this case as well.  \qed}
\end{mproof}
\vspace{-5pt}\noindent
We now are ready to prove that  

\vspace{-10pt}
$\left [Z_{S_{n}}\![2] \boxtimes Z_{S_{r}}(x_{1},x_2,\ldots,x_{r})\right ]\! (2,\!\ldots,\!2)\,\, \!\!\!\ge\! \left [Z_{S_{n}}[i] \boxtimes Z_{S_{r}}(x_{1},x_2,\ldots,x_{r})\right ]\! (2,\!\ldots,\!2), \forall i, 2\, \!\le\! i\!\le\!n!.$

\vspace{5pt}
\begin{theorem}
\label{SecondTermIsTheSecondBiggest}
{\rm 
$\,$\vspace{-8pt}
\begin{align}
& [Z_{S_{n}}[2] \boxtimes Z_{S_{r}}(x_{1},x_2,\ldots,x_{r})](2,2,\ldots,2) \geq [Z_{S_{n}}[i] \boxtimes Z_{S_{r}}(x_{1},x_2,\ldots,x_{r})](2,2,\ldots,2)\label{proposition5}
\end{align}

\vspace{-20pt}\noindent
$\forall i,2 \leq i \leq n!$ and $ \forall n, n < r$.
}
\end{theorem}

\begin{mproof}

{\rm  

Using Lemma~\ref{lemmaProp5} and Corollary~\ref{corollaryProp5} it suffices to show that
\vspace{-8pt}
\begin{align}
\label{prop5DerivedFormula}
& Z_{S_{r}}(2^{n-1},2^{n},\ldots) \geq Z_{S_{r}}(2^{\sum_{k=1}^{n}p_{i,k} {\rm gcd}(k,1)},\ldots,2^{\sum_{k=1}^{n}p_{i,k} {\rm gcd}(k,r)}).
\end{align}

\vspace{-15pt}\noindent
We prove the statement by induction on $r$.

\vspace{-10pt}\noindent
\textbf{Basis: $(r=1).$} By (\ref{cycleIndexRecurrence}), $Z_{S_{1}}(2^{n-1}) = 2^{n-1} Z_{S_{0}}() = 2^{n-1}.$ Similarly, by   (\ref{cycleIndexRecurrence}), $Z_{S_{1}}(2^{\sum_{k=1}^{n}p_{i,k} {\rm gcd}(k,1)}) = 
2^{\sum_{k=1}^{n}p_{i,k} {\rm gcd}(k,1)} Z_{S_{0}}() = 2^{\sum_{k=1}^{n}p_{i,k} }.$ Given that $\sum_{k=1}^{n}p_{i,k}\le n-1$ by Lemma~\ref{prop5Lemma2}, we have $2^{\sum_{k=1}^{n}p_{i,k} } \le 2^{n-1 },$ and hence the statement holds in this case.

\vspace{-5pt}\noindent
\textbf{Induction Step:} 
First, by (\ref{cycleIndexRecurrence}),
\begin{align}
Z_{S_{r}}(2^{n-1},2^{n},\ldots) &= 
\frac{1}{r} \begin{bmatrix}
2^{n-1} Z_{S_{r-1}}(2^{n-1},2^{n},\ldots) \\ 
 + 2^{n} Z_{S_{r-2}}(2^{n-1},2^{n},\ldots)\\ 
 + 2^{n-1} Z_{S_{r-3}}(2^{n-1},2^{n},\ldots)\\ 
 \vdots\\ 
 + 2^{\beta} Z_{S_{0}}()
\end{bmatrix}, 
\label{proposition5Eqna}
\end{align}

\vspace{-10pt}\noindent
where $\beta = n$ if $r$ is even and $\beta = n-1$ if $r$ is odd. Similarly,
\begin{align}
Z_{S_{r}}(2^{\sum_{k=1}^{n}p_{i,k} {\rm gcd}(k,1)},\ldots,2^{\sum_{k=1}^{n}p_{i,k} {\rm gcd}(k,r)}) =
\frac{1}{r} \begin{bmatrix}
2^{\sum_{k=1}^{n}p_{i,k} {\rm gcd}(k,1)} Z_{S_{r-1}}(2^{\sum_{k=1}^{n}p_{i,k} {\rm gcd}(k,1)},\ldots) \\
 + 2^{\sum_{k=1}^{n}p_{i,k} {\rm gcd}(k,2)} Z_{S_{r-2}}(2^{\sum_{k=1}^{n}p_{i,k} {\rm gcd}(k,1)},\ldots) \\
 + 2^{\sum_{k=1}^{n}p_{i,k} {\rm gcd}(k,3)} Z_{S_{r-3}}(2^{\sum_{k=1}^{n}p_{i,k} {\rm gcd}(k,1)},\ldots) \\
 \vdots\\ 
 + 2^{\sum_{k=1}^{n}p_{i,k} {\rm gcd}(k,r)} Z_{S_{0}}()
\end{bmatrix}
\label{proposition5Eqnb}
\end{align}

\vspace{-28pt}\noindent
Subtracting (\ref{proposition5Eqnb}) from (\ref{proposition5Eqna}), we have
\begin{align}
Z_{S_{r}}(2^{n-1},2^{n},\ldots) &- Z_{S_{r}}(2^{\sum_{k=1}^{n}p_{i,k} {\rm gcd}(k,1)},\ldots,2^{\sum_{k=1}^{n}p_{i,k} {\rm gcd}(k,r)})\nonumber \\=
\frac{1}{r} \begin{bmatrix}
2^{n-1} Z_{S_{r-1}}(2^{n-1},2^{n},\ldots) \\ 
 + 2^{n} Z_{S_{r-2}}(2^{n-1},2^{n},\ldots)\\ 
 + 2^{n-1} Z_{S_{r-3}}(2^{n-1},2^{n},\ldots)\\ 
 \vdots\\ 
 + 2^{\beta} Z_{S_{0}}()
\end{bmatrix} & -
\frac{1}{r} \begin{bmatrix}
2^{\sum_{k=1}^{n}p_{i,k} {\rm gcd}(k,1)} Z_{S_{r-1}}(2^{\sum_{k=1}^{n}p_{i,k} {\rm gcd}(k,1)},2^{\sum_{k=1}^{n}p_{i,k} {\rm gcd}(k,2)},\ldots) \\
 + 2^{\sum_{k=1}^{n}p_{i,k} {\rm gcd}(k,2)} Z_{S_{r-2}}(2^{\sum_{k=1}^{n}p_{i,k} {\rm gcd}(k,1)},2^{\sum_{k=1}^{n}p_{i,k} {\rm gcd}(k,2)},\ldots) \\
 + 2^{\sum_{k=1}^{n}p_{i,k} {\rm gcd}(k,3)}Z_{S_{r-3}}(2^{\sum_{k=1}^{n}p_{i,k} {\rm gcd}(k,1)},2^{\sum_{k=1}^{n}p_{i,k} {\rm gcd}(k,2)},\ldots) \\
 \vdots\\ 
 + 2^{\sum_{k=1}^{n}p_{i,k} {\rm gcd}(k,r)} Z_{S_{0}}()
\end{bmatrix}
\end{align}

\vspace{-25pt}\noindent
Thus, it suffices to show that the right hand side of the above equation is $\ge 0,$ or

\vspace{-25pt}
\begin{align}
\label{diffrenceEquation}
\begin{matrix}
2^{n-1} Z_{S_{r-1}}(2^{n-1},2^{n},\ldots) - 2^{\sum_{k=1}^{n}p_{i,k} {\rm gcd}(k,1)} Z_{S_{r-1}}(2^{\sum_{k=1}^{n}p_{i,k} {\rm gcd}(k,1)},2^{\sum_{k=1}^{n}p_{i,k} {\rm gcd}(k,2)},\ldots) \\ 
 + 2^{n} Z_{S_{r-2}}(2^{n-1},2^{n},\ldots) - 2^{\sum_{k=1}^{n}p_{i,k} {\rm gcd}(k,2)} Z_{S_{r-2}}(2^{\sum_{k=1}^{n}p_{i,k} {\rm gcd}(k,1)},2^{\sum_{k=1}^{n}p_{i,k} {\rm gcd}(k,2)},\ldots\\ 
 + 2^{n-1} Z_{S_{r-3}}(2^{n-1},2^{n},\ldots) - 2^{\sum_{k=1}^{n}p_{i,k} {\rm gcd}(k,3)} Z_{S_{r-3}}(2^{\sum_{k=1}^{n}p_{i,k} {\rm gcd}(k,1)},2^{\sum_{k=1}^{n}p_{i,k} {\rm gcd}(k,2)},\ldots)\\ 
 \vdots\\ 
 + 2^{\beta} Z_{S_{0}}() - 2^{\sum_{k=1}^{n}p_{i,k} {\rm gcd}(k,r)} Z_{S_{0}}()  \geq 0.
\end{matrix} &
\end{align}

\vspace{-10pt}\noindent Now  by induction hypothesis,  (\ref{prop5DerivedFormula}) holds for $1,2,\cdots, r-1.$ 
Thus, (\ref{diffrenceEquation}) can be replaced by

\vspace{-25pt}
\begin{align}
\label{diffrenceEquation2}
\begin{matrix}
2^{n-1} Z_{S_{r-1}}(2^{n-1},2^{n},\ldots) - 2^{\sum_{k=1}^{n}p_{i,k} {\rm gcd}(k,1)} Z_{S_{r-1}}(2^{n-1},2^{n},\ldots) \\ 
 + 2^{n} Z_{S_{r-2}}(2^{n-1},2^{n},\ldots) - 2^{\sum_{k=1}^{n}p_{i,k} {\rm gcd}(k,2)} Z_{S_{r-2}}(2^{n-1},2^{n},\ldots) \\ 
 + 2^{n-1} Z_{S_{r-3}}(2^{n-1},2^{n},\ldots) - 2^{\sum_{k=1}^{n}p_{i,k} {\rm gcd}(k,3)} Z_{S_{r-3}}(2^{n-1},2^{n},\ldots) \\ 
 \vdots\\ 
 + 2^{\beta} Z_{S_{0}}() - 2^{\sum_{k=1}^{n}p_{i,k} {\rm gcd}(k,r)} Z_{S_{0}}()  \geq 0.
\end{matrix} &
\end{align}

\vspace{-15pt}\noindent

\vspace{-15pt}\noindent
Moreover,  invoking Lemma~\ref{prop5Lemma2} gives

\vspace{-30pt}
\begin{align}
2^{n-1} Z_{S_{r-1}}(2^{n-1},2^{n},\ldots) & - 2^{\sum_{k=1}^{n}p_{i,k} {\rm gcd}(k,1)} Z_{S_{r-1}}(2^{n-1},2^{n}\ldots)\nonumber\\
& \geq 2^{n-1} Z_{S_{r-1}}(2^{n-1},2^{n},\ldots) - 2^{n-1} Z_{S_{r-1}}(2^{n-1},2^{n},\ldots) = 0. 
%& \geq 0
\end{align}

\vspace{-10pt}\noindent
Hence the difference in the first line in (\ref{diffrenceEquation2}) $\ge 0$, and 
therefore it is sufficient to show that

\vspace{-20pt}
\begin{align}
\label{inductionInequality10}
\begin{matrix}
2^{n} Z_{S_{r-2}}(2^{n-1},2^{n},\ldots) - 2^{\sum_{k=1}^{n}p_{i,k} {\rm gcd}(k,2)} Z_{S_{r-2}}(2^{n-1},2^n,\ldots) \\ 
 + 2^{n-1} Z_{S_{r-3}}(2^{n-1},2^{n},\ldots) - 2^{\sum_{k=1}^{n}p_{i,k} {\rm gcd}(k,3)} Z_{S_{r-3}}(2^{n-1},2^n,\ldots) \\ 
 \vdots\\ 
 %+ 2^{n-1} Z_{S_{1}}(2^{n-1}) - 2^{\sum_{k=1}^{n}p_{i,k} {\rm gcd}(k,r-1)} Z_{S_{1}}(2^{\sum_{k=1}^{n}p_{i,k} {\rm gcd}(k,1)})\\
  + 2^{\beta} Z_{S_{0}}() - 2^{\sum_{k=1}^{n}p_{i,k} {\rm gcd}(k,r)} Z_{S_{0}}()  \geq 0.
\end{matrix} &
\end{align}

\vspace{-18pt}\noindent To prove this inequality, we will combine four terms in pairs of consecutive lines for the remaining $r-1$ lines by considering two cases. If $r$ is odd then $\beta = n-1$ and no extra line remains in this pairing. Thus,  for all even $\alpha, 2 \leq \alpha \leq r-1,$  it suffices to prove

\vspace{-28pt}
\begin{align}
\label{inductionStep2}
\begin{matrix}
2^{n} Z_{S_{r-\alpha}}(2^{n-1},2^{n},\ldots) - 2^{\sum_{k=1}^{n}p_{i,k} {\rm gcd}(k,\alpha)} Z_{S_{r-\alpha}}(2^{n-1},2^n\ldots) ,\\ 
 + 2^{n-1} Z_{S_{r-\alpha-1}}(2^{n-1},2^{n},\ldots) - 2^{\sum_{k=1}^{n}p_{i,k} {\rm gcd}(k,\alpha+1)} Z_{S_{r-\alpha-1}}(2^{n-1},2^n\ldots)  \geq 0.
\end{matrix} &
\end{align}

\vspace{-20pt}\noindent
or,

\vspace{-30pt}
\begin{align}
\begin{matrix}
2^{n} Z_{S_{r-\alpha}}(2^{n-1},2^{n},\ldots) - 2^{\sum_{k=1}^{n}p_{i,k} k} Z_{S_{r-\alpha}}(2^{n-1},2^n\ldots) \\ 
 + 2^{n-1} Z_{S_{r-\alpha-1}}(2^{n-1},2^{n},\ldots) - 2^{\sum_{k=1}^{n}p_{i,k} {\rm gcd}(k,\alpha+1)} Z_{S_{r-\alpha-1}}(2^{n-1},2^n\ldots)  \geq 0.
\end{matrix} &
\end{align}

\vspace{-15pt}\noindent 
%Invoking Lemma~\ref{prop5Lemma3} gives
Now if  $\sum_{k=1}^{n}p_{i,k} {\rm gcd}(k,\alpha +1) \leq n-1 $, then 

\vspace{-20pt}
\begin{align}
%\begin{matrix}
%2^{n} Z_{S_{r-\alpha}}(2^{n-1},2^{n},\ldots) - 2^{\sum_{k=1}^{n}p_{i,k} {\rm gcd}(k,\alpha)} Z_{S_{r-\alpha}}(2^{n-1},2^{n},\ldots) \\ 
% + 2^{n-1} Z_{S_{r-\alpha-1}}(2^{n-1},2^{n},\ldots) - 2^{\sum_{k=1}^{n}p_{i,k} {\rm gcd}(k,\alpha+1)} Z_{S_{r-\alpha-1}}(2^{n-1},2^{n},\ldots)  \geq 0
%\end{matrix} & \\
\begin{matrix}
2^{n} Z_{S_{r-\alpha}}(2^{n-1},2^{n},\ldots) - 2^{\sum_{k=1}^{n}p_{i,k} k} Z_{S_{r-\alpha}}(2^{n-1},2^{n},\ldots) \\ 
 + 2^{n-1} Z_{S_{r-\alpha-1}}(2^{n-1},2^{n},\ldots) - 2^{n-1} Z_{S_{r-\alpha-1}}(2^{n-1},2^{n},\ldots)  \geq 
\end{matrix} &\nonumber \\
\begin{matrix}
2^{n} Z_{S_{r-\alpha}}(2^{n-1},2^{n},\ldots) - 2^{n} Z_{S_{r-\alpha}}(2^{n-1},2^{n},\ldots) \\ 
 + 2^{n-1} Z_{S_{r-\alpha-1}}(2^{n-1},2^{n},\ldots) - 2^{n-1} Z_{S_{r-\alpha-1}}(2^{n-1},2^{n},\ldots) = 0.
\end{matrix} &
\end{align}

\vspace{-12pt} \noindent 
On the other hand,  if $\sum_{k=1}^{n}p_{i,k} {\rm gcd}(k,\alpha+1) = n$, then we prove (\ref{inductionStep2}) by noting that  $\sum_{k=1}^{n}p_{i,k} {\rm gcd}(k,\alpha) \leq n-1$ by Lemma~\ref{prop5Lemma3}. Thus,

\vspace{-25pt}
\begin{align}
%\begin{bmatrix}
%2^{n} Z_{S_{r-\alpha}}(2^{n-1},2^{n},\ldots) - 2^{\sum_{k=1}^{n}p_{i,k} gcd(k,\alpha)} Z_{S_{r-\alpha}}(2^{n-1},2^{n},\ldots) \\ 
% + 2^{n-1} Z_{S_{r-\alpha-1}}(2^{n-1},2^{n},\ldots) - 2^{\sum_{k=1}^{n}p_{i,k} gcd(k,\alpha+1)} Z_{S_{r-\alpha-1}}(2^{n-1},2^{n},\ldots)
%\end{bmatrix} & \geq 0 \\
&2^{n} Z_{S_{r-\alpha}}(2^{n-1},2^{n},\ldots) - 2^{n-1} Z_{S_{r-\alpha}}(2^{n-1},2^{n},\ldots)\\
 &+  2^{n-1} Z_{S_{r-\alpha-1}}(2^{n-1},2^{n},\ldots) - 2^{n} Z_{S_{r-\alpha-1}}(2^{n-1},2^{n},\ldots)\nonumber\\
&= 2^{n-1} Z_{S_{r-\alpha}}(2^{n-1},2^{n},\ldots) - 2^{n-1} Z_{S_{r-\alpha-1}}(2^{n-1},2^{n},\ldots)\nonumber  \\
&2^{n-1} \left[Z_{S_{r-\alpha}}(2^{n-1},2^{n},\ldots) - Z_{S_{r-\alpha-1}}(2^{n-1},2^{n},\ldots)\right] 
%Z_{S_{r-\alpha}}(2^{n-1},2^{n},\ldots) - Z_{S_{r-\alpha-1}}(2^{n-1},2^{n},\ldots) & \geq 0
\end{align}

\vspace{-15pt}\noindent 
Now  by Lemma~\ref{prop5Lemma4}, $Z_{S_{r-\alpha}}(2^{n-1},2^{n},\ldots) \geq Z_{S_{r-\alpha-1}}(2^{n-1},2^{n},\ldots)$ and the statement is proved for odd $r, n < r.$ For even $r$, the last line in (\ref{inductionInequality10}) is left out in the pairing of consecutive lines and  $\beta = n$. In this case we have 
$2^{n} Z_{S_{0}}()  -2^{\sum_{k=1}^{n}p_{i,k} {\rm gcd}(k,r)} Z_{S_{0}}()  
 \geq 2^{n} Z_{S_{0}}() - 2^{\sum_{k=1}^{n}p_{i,k} k} Z_{S_{0}}() =
2^{n} Z_{S_{0}}() - 2^{n} Z_{S_{0}}()= 0$ and the statement follows. \qed
%\textbf{Case 2: r is an odd number} Odd case is shown to hold similarly and omitted. 

}

\end{mproof}

\vspace{10pt}
\begin{theorem}
\label{UpperBoundProp2}
\rm{
$\,$\\
\vspace{-25pt}

\begin{align}
& \left[ Z_{S_{n}}[2] \boxtimes Z_{S_{r}}(x_{1},x_2,\ldots,x_{r}) \right](2,2,\ldots,2) \leq \frac{\binom{r+2^{n}-1}{r}}{n! (n!-1)}.
\end{align}
}

\vspace{-30pt}\noindent
where $2 \leq n<r.$

\end{theorem}
\begin{mproof}
\rm{By Corollary~\ref{corollaryProp5}
\vspace{-10pt}
\begin{align}
& \left[Z_{S_{n}}[2] \boxtimes Z_{S_{r}}(x_{1},x_2,\ldots,x_{r})\right](2,2,\ldots,2) = \frac{1}{n!}Z_{S_{r}}(2^{n-1},2^{n},\ldots).
\end{align}

\vspace{-15pt}\noindent Thus, to prove the theorem, it is sufficient to show
\vspace{-12pt}
\begin{align}
\label{prop6InEq1}
\frac{1}{n!}Z_{S_{r}}(2^{n-1},2^{n},2^{n-1},2^{n},\ldots) \leq
\frac{\binom{r+2^{n}-1}{r}}{n! (n!-1)}
\end{align} 

\vspace{-20pt}\noindent 
where $2 \leq n<r$. 

\vspace{-5pt}\noindent
Now, using (\ref{cycleIndexRecurrence}), we get

\vspace{-25pt}
\begin{align}
\label{UpperBoundProp2Eq1}
& r Z_{S_{r}}(2^{n-1},2^{n},\ldots) =  \sum_{\textbf{odd } i}^{r-\beta_{1}}2^{n-1} Z_{S_{r-i}}(2^{n-1},2^{n},\ldots) + \sum_{\textbf{even } i}^{r-\beta_{2}}2^{n} Z_{S_{r-i}}(2^{n-1},2^{n},\ldots)
\end{align}

\vspace{-15pt}\noindent
where $\beta_{1}=1, \beta_{2}=0$ if $r$ is even and $\beta_{1}=0, \beta_{2}=1$ if $r$ is odd. Similarly, for $r-2$,

\vspace{-10pt}\noindent
\begin{align}
\label{UpperBoundProp2Eq2}
%& Z_{S_{r-2}}(2^{n-1},2^{n},\ldots) = \frac{1}{r-2} \left[ \sum_{\textbf{odd } i}^{r-2-\beta_{1}}2^{n-1} Z_{S_{r-2-i}}(2^{n-1},2^{n},\ldots) + \sum_{\textbf{even } i}^{r-2-\beta_{2}}2^{n} Z_{S_{r-2-i}}(2^{n-1},2^{n},\ldots) \right] \nonumber \\
&(r-2) Z_{S_{r-2}}(2^{n-1},2^{n},\ldots) = \sum_{\textbf{odd } i}^{r-2-\beta_{1}}2^{n-1} Z_{S_{r-2-i}}(2^{n-1},2^{n},\ldots) + \sum_{\textbf{even } i}^{r-2-\beta_{2}}2^{n} Z_{S_{r-2-i}}(2^{n-1},2^{n},\ldots).
\end{align}

\vspace{-25pt}\noindent
Subtracting (\ref{UpperBoundProp2Eq2}) from (\ref{UpperBoundProp2Eq1}) gives 

\vspace{-25pt}
\begin{align}
& r Z_{S_{r}}(2^{n-1},2^{n},\ldots) - (r-2) Z_{S_{r-2}}(2^{n-1},2^{n},\ldots)\nonumber \\
&= 2^{n-1} Z_{S_{r-1}}(2^{n-1},2^{n},\ldots) + 2^{n} Z_{S_{r-2}}(2^{n-1},2^{n},\ldots),  \\
& r Z_{S_{r}}(2^{n-1},2^{n},\ldots) = 2^{n-1} Z_{S_{r-1}}(2^{n-1},2^{n},\ldots) + (r-2+2^{n}) Z_{S_{r-2}}(2^{n-1},2^{n},\ldots), \\
& Z_{S_{r}}(2^{n-1},2^{n},\ldots) = \frac{1}{r} \left[2^{n-1} Z_{S_{r-1}}(2^{n-1},2^{n},\ldots) + (r-2+2^{n}) Z_{S_{r-2}}(2^{n-1},2^{n},\ldots)\right]. \label{recurrenceUpper1}
\end{align}

\vspace{-15pt} \noindent 
We will use  induction on $r$ and the recurrence given in (\ref{recurrenceUpper1}) to prove this inequality.

\vspace{-5pt}\noindent
\textbf{Basis. Case $r=3$:} Recall that

\vspace{-30pt}
\begin{align}
& Z_{S_{n}}[2] =  \frac{1}{n!} x_{1}^{n-2}x_{2}, \\
& Z_{S_{3}}(x_{1},x_{2},x_{3}) = \frac{1}{3!}(x_{1}^{3}+3x_{1}x_{2}+2x_{3}).
\end{align}

\vspace{-20pt}\noindent Thus,

\vspace{-30pt}
\begin{align}
& \left[ Z_{S_{n}}[2] \boxtimes Z_{S_{3}}(x_{1},x_{2},x_{3}) \right](2,2,\ldots,2)\nonumber \\
& = \left[ \frac{1}{n!} (x_{1}^{n-2}x_{2}) \boxtimes \frac{1}{3!}(x_{1}^{3}+3x_{1}x_{2}+2x_{3})\right](2,2,\ldots,2), \\
& = \frac{1}{3!n!} \left[(x_{1}^{n-2}x_{2}) \bigodot x_{1}^{3} + (x_{1}^{n-2}x_{2}) \bigodot (3x_{1}x_{2}) +(x_{1}^{n-2}x_{2}) \bigodot 2x_{3} \right](2,2,\ldots,2), \\
%& = \frac{1}{3!n!} \left[(x_{1}^{n-2}x_{2}) \bullet x_{1}^{3} + ((x_{1}^{n-2}x_{2}) \bullet 3x_{1}) ((x_{1}^{n-2}x_{2}) \bullet x_{2}) +(x_{1}^{n-2}x_{2}) \bullet 2x_{3} \right](2,2,\ldots,2), \\
& = \frac{1}{3!n!} \left[x_{1}^{3(n-2)}x_{2}^{3} + 3x_{1}^{n-2}x_{2}x_{2}^{n-2}x_{2}^{2} +
2x_{3}^{n-2}x_{6} \right](2,2,\ldots,2), \\
& = \frac{1}{3!n!} \left[2^{3n-3} + 3\times 2^{2n-1} + 2^{n} \right]\le \frac{\binom{r+2^{n}-1}{r}}{n! (n!-1)}.
\end{align}

\vspace{-20pt}\noindent
for $n = 2$ and $ r = 3.$

\vspace{-10pt}\noindent
\hskip 16pt {\bf Case} $r=4.$ In this case we have

\vspace{-30pt}
\begin{align}
& \left[ Z_{S_{n}}[2] \boxtimes Z_{S_{4}}(x_{1},x_{2},x_{3},x_{4}) \right](2,\ldots,2)\nonumber \\
& = \left[ \frac{1}{n!} (x_{1}^{n-2}x_{2}) \boxtimes \frac{1}{4!}(x_{1}^{4}+6x_{1}^{2}x_{2}+3x_{2}^{2}+
8x_{1}x_{3} + 6x_{4})\right](2,\ldots,2), \\
& = \frac{1}{4!n!} \bigg[(x_{1}^{n-2}x_{2}) \bigodot x_{1}^{4} + (x_{1}^{n-2}x_{2}) \bigodot (6x_{1}^{2}x_{2}) +(x_{1}^{n-2}x_{2}) \bigodot 3x_{2}^{2}\nonumber \\
& \quad +  (x_{1}^{n-2}x_{2}) \bigodot (8x_{1}x_{3}) +(x_{1}^{n-2}x_{2}) \bigodot 6x_{4} \bigg] (2,\ldots,2), \\
& = \frac{1}{4!n!} \left[x_{1}^{4(n-2)}x_{2}^{4} + 6x_{1}^{2(n-2)}x_{2}^{n-2}x_{2}^{2}x_{2}^{2} + 3x_{1}^{2(n-2)}x_{2}^{4} + 8x_{1}^{n-2}x_{3}^{n-2}x_{2}x_{6} + 6x_{4}^{n-2}x_{4}^{2} \right](2,\ldots,2), \\
& = \frac{1}{4!n!} \left[2^{4n-4} + 6 \times 2^{3n-2} + 3\times  2^{2n} + 8\times  2^{2n-2} + 6\times 2^{n} \right], \\
& = \frac{1}{4!n!} \left[2^{4n-4} + 6\times 2^{3n-2} + 5\times 2^{2n} + 6\times  2^{n} \right].
\end{align}

\vspace{-10pt}\noindent 
Now, given that $r = 4,$ the only possible values of $n$ are 2 and 3. If $n = 2$ then:

\vspace{-30pt}
\begin{align}
 \left[ Z_{S_{n}}[2] \boxtimes Z_{S_{4}}(x_{1},x_{2},x_{3},x_{4}) \right](2,2,\ldots,2) 
& = \frac{1}{4!n!} \left[2^{4n-4} + 6 \times 2^{3n-2} + 5\times  2^{2n} + 6 \times 2^{n} \right], \\
& = \frac{1}{4!2!} \left[2^{4} + 6 \times 2^{4} + 5\times  2^{4} + 6\times  2^{2} \right], \\
& = \frac{16+96+80+24}{4!2!} =  4.5,\\
&\le \frac{\binom{r+2^{n}-1}{r}}{n! (n!-1)} = \frac{\binom{7}{4}}{2! (2!-1)} =\frac{35}{2} = 17.5.
\end{align}

\vspace{-20pt}\noindent On the other hand, if $n = 3$ then:

\vspace{-30pt}
\begin{align}
 \left[ Z_{S_{n}}[2] \boxtimes Z_{S_{4}}(x_{1},x_{2},x_{3},x_{4}) \right](2,2,\ldots,2) 
& = \frac{1}{4!n!} \left[2^{4n-4} + 6\times 2^{3n-2} + 5 \times 2^{2n} + 6\times 2^{n} \right], 
\end{align}
\begin{align} 
& = \frac{1}{4!3!} \left[2^{8} + 6\times 2^{7} + 5 \times 2^{6} + 6\times 2^{3} \right],\\
& = \frac{256+768+320+48}{4!3!} = \frac{29}{3}, \\
& \le \frac{\binom{r+2^{n}-1}{r}}{n! (n!-1)} = \frac{\binom{11}{4}}{3! (3!-1)} =\frac{330}{30} = 11.
\end{align}

\vspace{-10pt}\noindent
\textbf{Induction Step:} Suppose that (\ref{prop6InEq1}) holds for all values from 3 to $r-1.$ Using the recurrence given in (\ref{recurrenceUpper1}) and the induction hypothesis for $r-1$ and $r-2$ we get:
\vspace{-5pt}
\begin{align}
 \frac{1}{n!}Z_{S_{r}}(2^{n-1},2^{n},\ldots) &= \frac{1}{n!r} \left[2^{n-1} Z_{S_{r-1}}(2^{n-1},2^{n},\ldots) + (r-2+2^{n}) Z_{S_{r-2}}(2^{n-1},2^{n},\ldots)\right], \\
& = \frac{2^{n-1}}{n!r} Z_{S_{r-1}}(2^{n-1},2^{n},\ldots) + \frac{r-2+2^{n}}{n!r}Z_{S_{r-2}}(2^{n-1},2^{n},\ldots),\\
& \leq \frac{2^{n-1}}{r}\frac{ \binom{r+2^{n}-2}{r-1}}{n! (n!-1)} + \frac{r-2+2^{n}}{r}\frac{\binom{r+2^{n}-3}{r-2}}{n! (n!-1)},
\end{align}
\begin{align}
& \leq \frac{2^{n-1}}  {n! (n!-1)r} \frac{(r+2^{n}-2)!}{(r-1)!(2^{n}-1)!} + \frac{r-2+2^{n}}{n! (n!-1)r} \frac{(r+2^{n}-3)!}{(r-2)!(2^{n}-1)!}, \\
& \leq \frac{2^{n-1}}  {n! (n!-1)r} \frac{(r+2^{n}-2)!}{(r-1)!(2^{n}-1)!} + \frac{(r-1) (r+2^{n}-2)!}{n! (n!-1)r!(2^{n}-1)!}, \\
& \frac{1}{n!}Z_{S_{r}}(2^{n-1},2^{n},\ldots) \leq \frac{(r+2^{n}-2)! (r+2^{n-1}-1)}{n! (n!-1)r!(2^{n}-1)!} \leq \frac{(r+2^{n}-2)! (r+2^{n}-1)}{n!(n!-1)r!(2^{n}-1)!}, \\
& \leq \frac{(r+2^{n}-1)!}{n! (n!-1)r!(2^{n}-1)!} = \frac{1}{n! (n!-1)} \binom{r+2^{n}-1}{r}, \\
& \leq \frac{1}{n! (n!-1)} \binom{r+2^{n}-1}{r}.
\end{align}

\vspace{-15pt}
This completes the proof. \qed
}
\end{mproof}

\vspace{-5pt}\noindent
Combining Theorems \ref{SecondTermIsTheSecondBiggest} and \ref{UpperBoundProp2}  concludes the upper bound calculation.

\vspace{5pt}
\begin{theorem}
\label{lastTheorem}
$|B_{u}(n,r)|\le  \frac{2\binom{r+2^{n}-1}{r}}{n!}.$
\end{theorem}
\begin{mproof}

$\,$\\
\vspace{-25pt}
\begin{align}
& \vert B_{u}(n,r)\vert = Z_{S_{n} \times S_{r}}(2,2,\ldots,2), \\
& = \left[Z_{S_{n}}(x_{1},x_2,\ldots,x_{n}) \boxtimes
Z_{S_{r}}(x_{1},x_2,\ldots,x_{r}) \right](2,2,\ldots,2),\\
& = \left[\left(Z_{S_{n}}[1] + Z_{S_{n}}[2] + \ldots + Z_{S_{n}}[n!]\right) \boxtimes
Z_{S_{r}}(x_{1},x_2,\ldots,x_{r})\right](2,2,\ldots,2), \\
& = \left[\left(Z_{S_{n}}[1] \right) \boxtimes
Z_{S_{r}}(x_{1},x_2,\ldots,x_{r})\right](2,2,\ldots,2) +  \left[\left(Z_{S_{n}}[2] \right) \boxtimes
Z_{S_{r}}(x_{1},x_2,\ldots,x_{r})\right](2,2,\ldots,2)\nonumber\\
&\quad +\ldots +\left[\left(Z_{S_{n}}[n!] \right) \boxtimes
Z_{S_{r}}(x_{1},x_2,\ldots,x_{r})\right](2,2,\ldots,2),\\
& \le \left[\left(Z_{S_{n}}[1] \right) \boxtimes
Z_{S_{r}}(x_{1},x_2,\ldots,x_{r})\right](2,2,\ldots,2) +  \left[\left(Z_{S_{n}}[2] \right) \boxtimes
Z_{S_{r}}(x_{1},x_2,\ldots,x_{r})\right](2,2,\ldots,2)\nonumber\\
&\quad +\ldots +\left[\left(Z_{S_{n}}[2] \right) \boxtimes
Z_{S_{r}}(x_{1},x_2,\ldots,x_{r})\right](2,2,\ldots,2), \\
& \leq \frac{\binom{r+2^{n}-1}{r}}{n!} + (n!-1) \frac{\binom{r+2^{n}-1}{r}}{n! (n!-1)} 
= \frac{2\binom{r+2^{n}-1}{r}}{n!}.\qed
\end{align}
\end{mproof}

\begin{remark}
{\rm 
It should be mentioned that, if $ r < n$,  using the relation $|B_u(n,r)| = |B_u(r,n)|$ gives
\begin{equation}
|B_u(n,r)|  \le  2\frac{\binom{n+2^{r}-1}{n}}{r!}.
\end{equation} 
Likewise,  if $ r < n,$ Theorem~\ref{theorem1-1} and $|B_u(n,r)| = |B_u(r,n)|$ together imply 
\begin{equation}
|B_u(n,r)|   \ge \frac{\binom{n+2^{r}-1}{n}}{r!}.
\end{equation} 

\vspace{-15pt}\noindent 
Furthermore, if $r = n,$ using the cycle index representation of bi-colored graphs provided in Section 3 in~\cite{harary1958number} and Theorem~\ref{theorem1-1} gives 
\begin{equation}
|B_u(n,n)|  \ge  \frac{\binom{n+2^{n}-1}{n}}{2n!}.
\end{equation} 

\vspace{-15pt}\noindent 
The $Z'$ term in the cycle index  representation of bi-colored graphs in \cite{harary1958number} prevents us from deriving an upper bound for $|B_u(n,n)|$ that is a constant multiple of the lower bound in this case. On the other hand, an obvious upper bound for $|B_u(n,n)|$ can be derived by setting $r = n+1$ in the inequality in Theorem~\ref{lastTheorem}.}
\end{remark}
\vspace{-10pt}\noindent

\vspace{-10pt}\noindent
{\bf Appendix:}

\vspace{-10pt}\noindent
 Table 1 lists $\ln|B_u(n,r)|$ along with the natural logarithms of lower and upper bounds for $1\le n < r \le 15.$

\vspace{5pt}
\bibliography{paperBib}{}
\bibliographystyle{unsrt}

\begin{table}[!t]
\centering{
\includegraphics[scale=0.4]{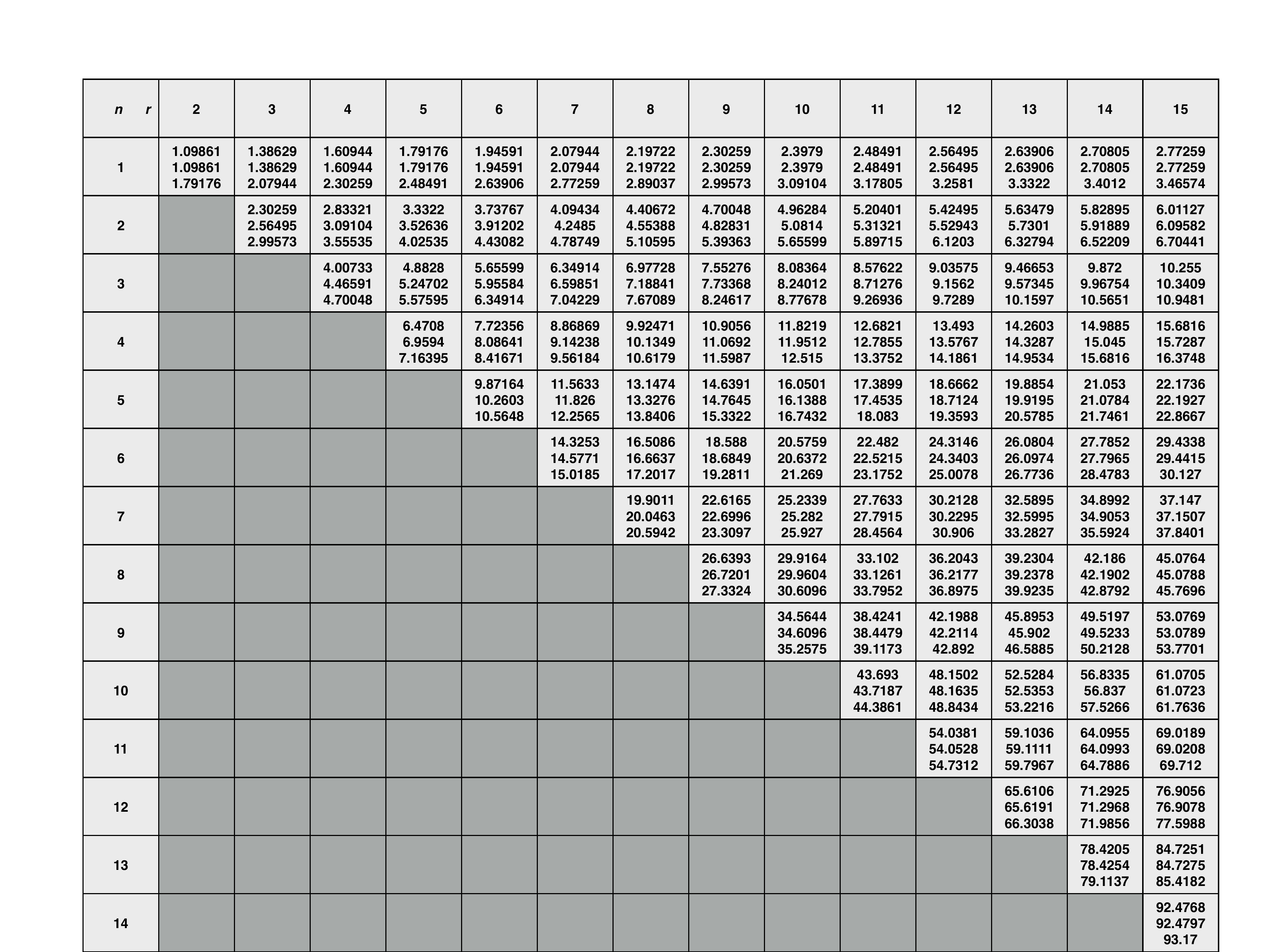}
\caption{Exact values of $\ln |B_u(n,r)|,1\le n < r \le 15,$  and natural logarithms of lower and upper bounds.}
\label{boundsAndExactValues}
}
\end{table} 

\end{document}